\documentclass[a4paper]{article}

\usepackage{amsmath}
\usepackage{amssymb}
\usepackage{amsthm}
\usepackage[cp1250]{inputenc}
\usepackage[T1]{fontenc}
\usepackage[english]{babel}
\usepackage{graphicx}
\usepackage{amsfonts}

 \newtheorem{thm}{Theorem}[section]
 \newtheorem{cor}[thm]{Corollary}
 \newtheorem{lem}[thm]{Lemma}
 \newtheorem{prop}[thm]{Proposition}
 \theoremstyle{definition}
 
 \theoremstyle{remark}
 \newtheorem{rem}[thm]{Remark}
 
 \numberwithin{equation}{section}

\begin{document}

\title{Asymptotic behavior of solutions of the Dirac system
  with an integrable potential}
  \date{August 26, 2020}
\maketitle

\author{\textbf{Łukasz Rzepnicki}\\
Faculty of Mathematics and Computer Science\\
Nicolaus Copernicus University\\
ul. Chopina 12/18, 87-100 Toruń\\
Poland\\ \verb"keleb@mat.umk.pl"}

\begin{abstract}
We consider the  Dirac system on the interval $[0,1]$ with a spectral parameter $\mu\in\mathbb{C}$ and
a complex-valued potential with entries from $L_p[0,1]$, where $1\leq p <2$.
We study the asymptotic behavior of its solutions  in a stripe $|{\rm Im}\,\mu|\le d$ for $\mu\to \infty$.
These results allows us to obtain  sharp
asymptotic formulas for eigenvalues and eigenfunctions
of  Sturm--Liouville operators associated with the aforementioned Dirac system.
\end{abstract}

\textbf{keywords:}\\
Dirac system, spectral problem, integrable potential,
Sturm--Liouville operator\\
\verb"MSC"[2010] Primary 34L20,  Secondary 34E05

\section*{Acknowledgement}
The author was supported by NCN grant no. UMO-2017/27/B/ST1/00078.

\section{Introduction}

Let consider for $x\in [0,1],$ a Cauchy
problem
\begin{equation}\label{v120}
D'(x)+J(x) D(x)=A_\mu D(x),\quad
D(0)=I,
\end{equation}
where  $A_\mu=i\mu J_0$,
and
\begin{equation}\label{matrixA}
J_0=\left[\begin{array}{cc}
1 & 0\\
0 &-1
\end{array}
\right],\quad
J(x)=\left[\begin{array}{cc}
0 & \sigma_1(x)\\
\sigma_2(x) & 0
\end{array}
\right],\quad
I:=\left[\begin{array}{cc}
1 & 0 \\
0 & 1
\end{array}\right],
\end{equation}
$\mu\in \mathbb{C}$ is a spectral parameter, and for $j=1,2$ complex-valued functions $\sigma_j$ belong to  $L_p[0,1]$, where $1\leq p < 2$.
We study the asymptotic behavior of its solutions $D(x)=D(x,\mu)$ with respect to $\mu$ from a horizontal stripe
 \[
P_d:=\{\mu\in \mathbb{C}:\,|{\rm Im}\,\mu|\le d\}.
\]
and $\mu\to \infty$.

The solution of \eqref{v120}, is  a
matrix $D$ with entries from
the space of absolutely continuous on $[0,1]$ functions
(i.e. from the Sobolev space $W_1^1[0,1]$)
satisfying  \eqref{v120} for a.e. $x\in [0,1]$.
In our case, this conditions together with the equation yield that $D$ has entries from $W^1_p[0,1]$.

This article is an addendum to the paper \cite{GRZ1}, where the problem
\eqref{v120} was analyzed for $\sigma_j\in L_2[0,1]$, $j=1,2$. In that text one can find
 background for Dirac systems and their connection with Sturm--Liouville problems.

We relay here on the same method as in \cite{GRZ1} to use all of its  advantages and obtain sharp
asymptotic formulas for $D$ and consequently for  spectral problems associated with \eqref{v120}.
In the case when $\sigma_j\in L_p[0,1]$, $j=1,2$, $ p > 2$ one can use the results from \cite{GRZ1} due to the obvious embedding
between $L_p[0,1]$ spaces.
Thus, in this text we restrict ourselves only to $1\leq p < 2$.

We are interested in the following spectral problem:
\begin{equation}\label{Dop}
Y'(x)+J(x) Y(x)=A_\mu Y(x),\quad
 x\in [0,1],
\end{equation}
where $Y=[y_1,y_2]^T$ and
\begin{align}\label{bound12}
y_1(0)=y_2(0), \ \ y_1(1)=y_2(1).
\end{align}

 Conditions \eqref{bound12}  are
 an example of \textit{strongly regular} boundary conditions. The Dirac-type systems or equation \eqref{Dop} with a general formulation
of \textit{regular}  or \textit{strongly regular} conditions have been studied recently in
many papers and different method.

In \cite{SSDirac} A. M. Savchuk and A. A. Shkalikov  derived for $p\geq1$ basic asymptotic formulas for eigenvalues and for fundamental solutions
  of the Dirac-type system only with the leading term and the reminders expressed by
  $\gamma_0$ and $\gamma$ given by \eqref{ga0} and \eqref{ga10}. They obtained their
 results applying Pr\"{u}fer’s substitution.

  Their result is equivalent to  first thesis \eqref{t222} of corollary \ref{Cor3.2}.
  Note that next statement \eqref{LCdop} is a significant extension of the previous result.
   Its version for $p=1$ may be found in remark \ref{p1}. The most general result is the content of lemma \ref{asD}.

  Using our method it is also possible to obtain very detailed
  formulas for eigenvalues and eigenfunctions. In case of  the spectral problem associated with \eqref{bound12} the eigenvalues
  admit the representation \eqref{mu01}-\eqref{mu0} with remainders satisfying \eqref{forp1} and \eqref{forp} for $p=1$ and
  $1<p<2$ respectively. In literature (for instance in \cite{SSDirac}) for $1<p<2$ one may found results which state that eigenvalues are of the form
  $\pi n + r_n$, where $(r_n)\in l_q$, and $q$ is conjugated to $p$.
  Here it is worth to underline that  beside the leading term in our asymptotic formulas there
  occur Fourier coefficients of known functions and the reminder, which belongs to $l_{q/2}$.
  Additionally, for $p=1$ we extend known formulas with  $|r_n|<c \Gamma(\pi n)$ (where
  $\Gamma$ is defined in \eqref{Gam}) into more detailed one  with the remainder satisfying $|r_n|<c \Gamma^2(\pi n)$.

  In the same spirit theorem \ref{more} and corollary \ref{less} related to eigenfunctions generalize significantly those from literature.

  Our method is applicable not only to the spectral problem  \eqref{Dop}-\eqref{bound12} but it works
  as well for different cases of strongly regular boundary conditions. What is more it may be used
 to deal with the class of \textit{regular} boundary conditions (in the sense of Birkhoff).

   The articles of A. M. Savchuk and I. V. Sadovnichaya: \cite{Sad2}, \cite{S}, \cite{SSad2} and \cite{SSad1}
 may be regarded as a continuation of method from \cite{SSDirac} and its application for $p=1$ to problems from the fields of asymptotics formulas and basis properties.
 Almost all aforementioned works prove or use the same type of results as mentioned before since they deal with the Riesz basis property and very detailed formulas are not needed.

 The same aims had M. M. Malamud, A. V. Agibalova and L. L. Oridoroga in \cite{AMO} for $p=2$ and latter two in \cite{LMal} for $p=1$. Here the authors
used the method of transformation operators.

Whereas in order to study  inverse spectral problems S. Albeverio, R. Hryniv and Y. Mykytyuk in \cite{Alb1} investigated a
direct spectral problem for the Dirac system in the form
\begin{equation}\label{System3}
BZ'(x)+Q(x)Z(x)
=\mu Z(x),\quad x\in [0,1],
\end{equation}
where
\[
B=\left[\begin{array}{cc}
0 & 1\\
-1 &  0
\end{array}
\right],\quad  Q(x)=\left[
\begin{array}{cc}
q_{1}(x) & q_{2}(x)\\
q_{2}(x) & -q_{1}(x)\end{array}
\right],\quad q_{j}\in L_p[0,1], \ \ j=1,2,
\]
with $p\geq1$. They proved also short formulas for fundamental system of solutions, where
reminders were expressed in terms of Fourier coefficients for unknown functions from $L_p$. Furthermore, for the operators associated with
 the system \eqref{System3} with two kinds of conditions
\begin{align}\label{boundij}
z_j(1)=z_2(0)=0, \ \ j=1,2,
\end{align}
they presented basic formulas for eigenvalues with the same type of reminders.
That class of results can be directly derived from our approach with the help of
transformation $Z=UY$, where
$$U=\left[\begin{array}{cc}
1 & -i\\
-i &  1
\end{array}\right].
$$
It leads to the
system \eqref{v120} with $\sigma_1=q_1+iq_2$
and $\sigma_2=q_1-iq_2$ with appropriate conditions.
The relation between different formulations of Dirac systems is explained deeper
in \cite{GRZ1}.

More results concerning different type of problems for the Dirac system with  may be found in the series of paper of
P. Djakov and B. Mityagin: \cite{DM0}, \cite{DM2} and \cite{DM} or D. V. Puyda \cite{Puyy}.

We start with the section concerning asymptotic behavior for solutions of Dirac system. Next, in section 3 we apply these results to the aforementioned spectral problem. For the clarity of exposition some technical results are placed at the end in appendix.

\section{Dirac system and its solutions}

In this section we study the matrix Cauchy problem \eqref{v120} and the behavior of its solution in a special integral form.
The idea of this approach was taken from
\cite[Ch. 1, $\S 24$]{March} and developed in \cite{GRZ1}.
We follow it here directly for similar operators but in different function spaces.

First, we introduce a necessary notation. We use throughout the text a standard symbol
 $L_p[0,1]$, $p\geq1$ to denote the space of measurable complex functions integrable with $p$-th power with the classical norm
$$\|f\|_{L_p}=\Big(\int_{0}^1|f(x)|^p dx\Big)^{1/p}.$$
We write $l_p$, $p\geq1$ for the space of complex sequences summable with
$p$-th power and endowed with the norm $$\|(x_n)\|_p=\Big(\sum_{n=1}^{\infty}|x_n|^p\Big)^{1/p}.$$ $W_p^1[0,1]$ is a standard Sobolev space with
the derivative in $L_p[0,1]$.

If $X$ is a Banach space, then  $M(X)$ stands for  the Banach
space
of $2\times 2$ matrices with entries
from $X$ and the norm
\[
\|Q\|_{M(B)}:=\sum_{k,j=1}^2 \|Q_{jk}\|_B,\quad
Q=[Q_{jk}]_{j,k=1}^2.
\]

We assume throughout the text that $1\leq p <2$. Morever if $1<p<2$, then let $q$ and $p$ be conjugate exponents and $r$
 be the number from Young's  convolution  inequality  ie.
\begin{align}\label{qr}
\frac{1}{p}+\frac{1}{q}=1\ \ \mbox{and} \ \  r=\frac{2}{2-p}.
\end{align}
If $p=1$, then $q=\infty$ and $r=1$.
Let
\begin{align}\label{Delta}
\Delta:=\{(x,t)\in\mathbb{R}^2: \; 0\leq t\leq x \leq 1\}
\end{align}
and
\begin{align*}
B:=\{f\colon [0,1]\times [0,1]\to \mathbb{C}\; a.e.\; &:\forall_{x\in [0,1]} \; f(x,\cdot)\in C([0,1],L_r),\\ &\mathrm{supp} f\subset \Delta)\}.
\end{align*}
We equip $B$
with the norm
$$\|f\|_B:=\sup_{x\in [0,1]}\|f(x,\cdot)\|_{L_r[0,x]},$$
so that $B$ is a Banach space. In particular, directly from the definition if $f\in B$, then $f(x,t)=0$ for $0\leq x < t \leq 1$.
This comment allows us  to underline the property which will be used in the text i.e.
for $f\in B$  there holds
\begin{align}\label{zerowanie}
\int_0^xf(x,t)dt=\int_0^1f(x,t)dt \in C[0,1].
\end{align}

We will use the series of constants connected with functions $\sigma_j$, $j=1,2$ in our estimations:
\begin{equation}\label{coetants}
a_0:=\max\{\|\sigma_1\|_{L_1},\|\sigma_2\|_{L_1}\},\quad
a:=\|\sigma_1\|_{L_1}\cdot\|\sigma_2\|_{L_1},\quad
a_1:=\|\sigma_1\|_{L_1}+\|\sigma_2\|_{L_1},
\end{equation}
and
\begin{equation}\label{constants1}
\tilde{a}_0:=\max\{\|\sigma_1\|_{L_p},\|\sigma_2\|_{L_p}\},\quad
\tilde{a}:=\|\sigma_1\|_{L_p}\cdot\|\sigma_2\|_{L_p},
\quad
a_2:=\|\sigma_1\|_{L_p}+\|\sigma_2\|_{L_p}.
\end{equation}
Moreover, let
\begin{equation}\label{NfunC}
\sigma_0(x):=|\sigma_1(x)|+|\sigma_2(x)|\in L_p[0,1].
\end{equation}

Now we are ready to establish  a first crucial property of the solutions to \eqref{v120}. The proof of the following lemma relays on
technical results
related to certain integral operators, which are placed in appendix.
\begin{lem}\label{jednor}
Let $\sigma\in L_p[0,1]$, $1\leq p < 2$.
\begin{itemize}
\item[a)]  The unique solution $D=D(x,\mu)$ of Cauchy problem \eqref{v120} 
can be represented as
\begin{align}\label{D}
D(x,\mu)=e^{xA_\mu}+\int_0^xe^{(x-2t)A_\mu}[J(t)+Q(x,t)]dt,
\end{align}
where
$Q\in M(B)$ is the unique solution of the integral equation
\begin{align}\label{Qtilde}
Q(x,t)=
\tilde{J}(x,t)+
\int_0^{x-t}J(t+\xi)Q(t+\xi,\xi)d\xi,
\end{align}
with $\tilde{J}\in M(B)$  is given by
\begin{equation}\label{ksigmazeroo}
\tilde{J}(x,t):=
\int_0^{x-t}J(t+\xi)J(\xi)d\xi=
\int_t^x J(s)J(s-t)\,ds,\quad(x,t)\in \Delta.
\end{equation}
\item[b)] The following estimates hold:
\begin{equation}\label{L1Est}
\|Q\|_{M(B)}\le
c, \ \ \|D\|_{M(C[0,1])}
\leq c, \ \ \mu\in P_d
\end{equation}
with certain constants $c=c(d,\sigma_1,\sigma_2)$.
\end{itemize}
\end{lem}

\begin{proof}
Note that the uniqueness of solutions comes from  general results on Sturm--Liouville equations (for instance \cite[Thm. 1.2.1]{Zettl}). We look for  solutions
of \eqref{v120} in a special form
\begin{align}\label{wzor1}
D(x,\mu)=e^{xA_\mu}U(x,\mu),
\quad U(0,\mu)=I.
\end{align}
The identity
\begin{equation}\label{relAA}
J(x)e^{xA_\mu}=e^{-xA_\mu}J(x),\quad \mbox{a. e.}\;\;x\in [0,1]
\end{equation}
yield that $U$ satisfies the Cauchy problem
\[
U'(x,\mu)+e^{-2x A_\mu}J(x) U(x,\mu)=0,\qquad x\in [0,1],\quad U(0,\mu)=I,
\]
and this is equivalent to the integral equation
\begin{align}\label{wzor15}
U(x,\mu)=I-\int_0^xe^{-2tA_\mu}J(t)U(t,\mu)dt,\qquad x\in [0,1].
\end{align}
We will seek for solutions of
\eqref{wzor15}
in the form
\begin{align}\label{wzor2}
U(x,\mu)=I+\int_0^xe^{-2tA_\mu}Q_0(x,t)dt,\quad
\end{align}
 where $Q_0\in M(B)$ does not depend on $\mu$.
 Inserting \eqref{wzor2} into \eqref{wzor15},
we obtain
\begin{align*}
\int_0^x e^{-2tA_\mu}Q_0(x,t)\,dt&=
-\int_0^x e^{-2tA_\mu}J(t)\,dt\\
&-\int_0^x e^{-2tA_\mu}J(t)
\int_0^t e^{-2s A_\mu} Q_0(t,s)\,ds\,dt.
\end{align*}
Due to the fact
\begin{equation}\label{Note1}
J_0^2=I,\quad J_0 J(x)+J(x)J_0=0,\qquad \mbox{a.e.}\;x\in [0,1],
\end{equation}
we get
\begin{align*}
\int_0^xe^{-2tA_\mu}J(t)
\int_0^t
e^{-2s A_\mu} &Q_0(t,s)ds\,dt
=\int_0^xe^{-2tA_\mu}
\int_0^t
e^{2s A_\mu} J(t)Q_0(t,s)ds\,dt
\\
&=\int_0^xe^{-2tA_\mu}\int_0^{x-t}
J(t+\xi)Q_0(t+\xi,\xi)d\xi dt,
\end{align*}
thus
\[
\int_0^xe^{-2tA_\mu}Q_0(x,t)\,dt=
-\int_0^xe^{-2tA_\mu}\left(J(t)
+\int_0^{x-t}
J(t+\xi) Q_0(t+\xi,\xi)d\xi\right)\,dt
\]
for all $x\in [0,1].$ We conclude that $U$ is a solution of \eqref{wzor15}
if and only if $Q_0\in M(B)$ is a solution of
\begin{align}\label{prev}
Q_0(x,t)=-J(t)-\int_0^{x-t}J(t+\xi)Q_0(t+\xi,\xi)d\xi.
\end{align}
Next,  setting
\[
Q_0(x,t)=-J(t)+Q(x,t),\qquad (x,t)\in \Delta,
\]
and using \eqref{Qtilde}, we infer that  $Q$ satisfies \eqref{prev}.
For $Q$ the equation \eqref{prev} can be rewritten in an operator form
\[
Q=\tilde{J}+\tilde{T} Q,\quad
\tilde{T}=
-\left[\begin{array}{cc}
0 & T_{\sigma_1}\\
T_{\sigma_2} & 0
\end{array}\right],
\]
for the operators $T_{\sigma_1}$ and $T_{\sigma_2},$ defined on $B$ by
\begin{align}\label{T}
(T_{\sigma}f)(x,t)=\int_0^{x-t}\sigma(t+\xi)f(t+\xi,\xi)d\xi=
\int_t^x\sigma(s)f(s,s-t)ds, 
\end{align}
where $\sigma\in L_p[0,1]$.

Observe that
\[
\tilde{J}(x,t)=\left(\begin{array}{cc}
\tilde{\sigma}_1(x,t) & 0\\
0 & \tilde{\sigma}_2(x,t)\end{array}\right),
\]
where
\begin{equation}\label{defSS}
\tilde{\sigma}_1(x,t):=\int_0^{x-t}\sigma_1(t+\xi)\sigma_2(\xi)d\xi,\quad
\tilde{\sigma}_2(x,t):=\int_0^{x-t}\sigma_2(t+\xi)\sigma_1(\xi)d\xi.
\end{equation}
 According to lemma \ref{Jtilde} $\tilde{J}\in M(B)$. What is more, the operators $T_{\sigma}$ are  linear and bounded  on $B$ due to lemma \ref{Tsigma}.
In particular,  we have
\[
\|\tilde{T}F\|_{M(B)}
\le
a_0
\|F\|_{M(B)},\quad F\in M(B).
\]

Next observe that
\[
\tilde{T}^{2n}=
\left[\begin{array}{cc}
T_{12}^n & 0\\
0 & T_{21}^n
\end{array}\right],\quad n\in {\rm N},
\]
for bounded linear operators
$T_{12}$ and $T_{21}$ on $B$ given by
\[
T_{12}:=T_{\sigma_1} T_{\sigma_2},\quad
T_{21}:=T_{\sigma_2} T_{\sigma_1}.
\]
Therefore by \eqref{Tn}, we derive
\[
\|\tilde{T}^{2n}F\|_{M(B)}
\le
\frac{a^n}{n!}\|F\|_{M(B)},\quad F\in M(B).
\]
We thus  see that
\eqref{Qtilde} has a unique solution
$Q\in M(B)$ of the form
\begin{equation}\label{Pred}
Q=\sum_{n=0}^\infty \tilde{T}^n
\tilde{J}=\sum_{n=0}^\infty \tilde{T}^{2n}(I+\tilde{T})\tilde{J},
\end{equation}
and moreover
\begin{equation}\label{Pred11}
\quad \|Q\|_{M(B)}
\le (1+a_0)e^a\|\tilde{J}\|_{M(B)}.
\end{equation}
Then \eqref{Pred} and \eqref{L22CY} imply \eqref{L1Est}.

Note that from \eqref{D} and \eqref{zerowanie} we have $D\in C[0,1]$.
Adding together of \eqref{D} and \eqref{L1Est}, we obtain
\begin{equation}\label{LLl1}
\|D\|_{M(C[0,1])}\leq e^{d} \Big(1+a_1+\|Q(x,t)\|_{M(B)}\big),\quad \mu\in P_d.
\end{equation}
\end{proof}

We now proceed to derivation of asymptotic formulas for $D$ with the use of the previous lemma.
In what follows we will use different types of
estimates for reminders.
For fixed $\sigma_j\in L_p$, $p\geq1$, $j=1,2$, and  $\mu\in \mathbb{C}$ define
\begin{align}\label{ga0}
\gamma(\mu)&:=\sum_{j=1}^2\left(\Big\|\int_0^xe^{-2i\mu t}\sigma_j(t)dt\Big\|_{L_q}+\Big\|\int_0^xe^{2i\mu t}\sigma_j(t)dt\Big\|_{L_q}
\right),
\end{align}
where $1/q+1/p=1$.
We will need also
\begin{align}\label{ga10}
\gamma_0(x,\mu)&:=\sum_{j=1}^2\left(\Big|\int_0^xe^{-2i\mu t}\sigma_j(t)dt\Big|+\Big|\int_0^xe^{2i\mu t}\sigma_j(t)dt\Big|\right),\quad x\in [0,1]
\end{align}
and
\begin{equation}\label{gamma12}
\gamma_1(\mu):=\int_0^1 \sigma_0(s)\gamma_0^2(s,\mu)\,ds,\qquad
\gamma_2(\mu):=l^2_2\gamma^2(\mu)+l_1\gamma_1(\mu).
\end{equation}
\begin{align}\label{Gam}
\Gamma(\mu)&:=\sum_{j=1}^2\left(\sup_{x\in[0,1]}\Big|\int_0^xe^{-2i\mu t}\sigma_j(t)dt\Big|+\sup_{x\in[0,1]}\Big|\int_0^xe^{2i\mu t}\sigma_j(t)dt\Big|\right).
\end{align}
Note that $\Gamma$ is nothing else than $\gamma$ for $p=1$ and $q=\infty$.

It is easy to see that if $\mu\in P_d$ then
\begin{equation}\label{gS1}
\gamma_0(x,\mu)\le 2e^{2d}a_1, \ \ \|\gamma_0(x,\mu)\|_{L_q}\le \gamma(\mu),\quad
\gamma(\mu)\le 2e^{2d}a_1,
\quad x\in [0,1],
\end{equation}
and
\begin{equation}\label{gamma2-1}
\gamma_2(\mu)\le
4e^{4d}a_1^2(a_2+a_1^2),\quad \gamma_2(\mu)\le 2a_1e^{2d}(a_2+2a_1e^{2d}\|\sigma_0\|_{L_p})\gamma(\mu),
\end{equation}

In the following lemma we will need
\begin{equation}\label{NotNN}
N(x,t):=(\tilde{J}+\tilde{T}\tilde{J})(x,t)\in B.
\end{equation}
Observe  that the explicit form of $N$ is
\[
N(x,t)= \left(\begin{array}{cc}
\tilde{\sigma}_1(x,t) & -(T_{\sigma_1}\tilde{\sigma}_2)(x,t)\\
-(T_{\sigma_2}\tilde{\sigma}_1)(x,t) &  \tilde{\sigma}_2(x,t)
\end{array}
\right).
\]

The very basic but crucial result use mainly the description of some integrals connected with the operator $\tilde{T}$
and its powers stated in lemma \ref{Asimp}.
\begin{lem}\label{asD}
Let  $\sigma_j\in L_p$, $1\leq p <2$ for $j=1,2$.
If $D(x,\mu)$ is a solution of \eqref{v120} then
\begin{align}\label{t1}
D(x,\mu)&=e^{x A_\mu} + D^{(0)}(x,\mu)+D^{(1)}(x,\mu),
\end{align}
where
\[
D^{(0)}(x,\mu)=\int_0^x e^{(x-2t)A_\mu} J(t)\,dt
+\int_0^x
e^{(x-2t)A_\mu} N(x,t)\,dt,
\]
and for all $\mu\in P_d$ and $x \in [0,1],$
\[
\|D^{(1)}(x,\mu)\|_{M(C[0,1])}\le c\gamma_2(\mu),
\]
where $c=c(d,\sigma_1,\sigma_2)$.
\end{lem}

\begin{proof}
Going back to the formulas \eqref{D} and
\eqref{Pred} for
$D=D(x,\mu)$, $x\in [0,1]$, $\mu\in P_d$,
note that \begin{align}
D(x,\mu)&=e^{xA_\mu}+\int_0^xe^{(x-2t)A_\mu}[J(t)+Q(x,t)]dt
\nonumber\\
&=e^{xA_\mu}+\int_0^xe^{(x-2t)A_\mu}J(t)dt
+\int_0^xe^{(x-2t)A_\mu}\tilde{J}(x,t)dt
\nonumber\\
&+\int_0^xe^{(x-2t)A_\mu}(\tilde{T}\tilde{J})(x,t)dt+D^{(1)}(x,\mu)\label{DDD}
\end{align}
where
\[
D^{(1)}(x,\mu)=\int_0^xe^{(x-2t)A_\mu}\sum_{n=2}^\infty (\tilde{T}^n\tilde{J})(x,t)dt.
\]
Using \eqref{DDD} and the inequality \eqref{EstIm3} proved in appendix, we infer that
\begin{align*}
\|D^{(1)}(x,\mu)\|_{M(C[0,1])}
&\le \sum_{n=2}^\infty \left\|\int_0^xe^{(x-2t)A_\mu}(\tilde{T}^n\tilde{J})(x,t)dt\right\|_{M(C[0,1])}
\\
&\le 2 e^d \gamma_2(\mu)\sum_{n=2}^\infty \frac{e^{2nd}a_1^{n-2}}{(n-2)!}=
2e^{5d}\exp{(e^{2d}a_1)}\gamma_2(\mu),
\end{align*}
for all $x\in [0,1]$ and $\mu\in P_d$.
\end{proof}

The above lemma leads to sharp asymptotic formulas for $D$, which are the main result of this section.
\begin{cor}\label{Cor3.2}
For every  $d>0$ there exist  $c_j=c_j(d,\sigma_1,\sigma_2)$, $j=0,1,$ such that
for all $x\in [0,1]$ and $\mu\in P_d,$
\begin{equation}\label{t222}
D(x,\mu)=e^{x A_\mu}+R(x,\mu),
\end{equation}
where
\[
\|R(x,\mu)\|_{M(C[0,1])}\le c_0,\quad
\|R(x,\mu)\|_{M(\rm C)}\le
 c_1(\gamma(\mu)+\gamma_0(x,\mu)).
\]
Moreover,
\begin{equation}\label{LCdop}
D(x,\mu)=e^{xA_\mu}
+D_0(x,\mu)+R_0(x,\mu),
\end{equation}
where
\[
D_{0}(x,\mu):=\int_0^x e^{(x-2t)A_\mu}(J(t)+\tilde{J}(t))\,dt
\]
and
\[
\|R_0(x,\mu)\|_{M(\mathbb{C})}\le
 c_2(\gamma(\mu)\gamma_0(x,\mu)+\gamma_2(\mu)),\quad x\in [0,1].
\]
\end{cor}
\begin{proof}
Let us start with  several
simple observations.
First of all, remark that
\begin{align*}
\|D^{(0)}(x,\mu)\|_{M(C([0,1]))}&\le e^d(\|J\|_{M(L_1[0,1])}
+\|\tilde{J}\|_{M(C(\Delta))}
+\|\tilde{T}\tilde{J}\|_{M(C(\Delta))})
\\
&\le e^d(a_1+(1+a_1)\|\tilde{J}\|_{M(C(\Delta))})\\
&\le
e^d(a_1+2(1+a_1)\tilde{l}).
\end{align*}
Note also that from
\[
\left\|\int_0^x e^{(x-2t)A_\mu} J(t)\,dt\right\|_{M(\mathbb{C})}\le
e^d\gamma_0(x,\mu),\quad x\in [0,1],
\]
\eqref{EstIm0}, and \eqref{EstIm100} it follows that
\[
\left\|D^{(0)}(x,\mu)\right\|_{M({\rm C}}\le
e^d\gamma_0(x,\mu)+ 2e^{5d} (1+\tilde{a}_0)\tilde{a}_0 \gamma(\mu),\quad x\in [0,1].
\]
Furthermore, by \eqref{EstIm1},
\begin{align*}
\left\|\int_0^x e^{i\mu(x-2t)} \tilde{T}\tilde{J}(x,t)\,dt\right\|_{M(\mathbb{C})}
&\le
e^d\left\|\int_0^x e^{-2i\mu t} \tilde{T}\tilde{J}(x,t)\,dt\right\|_{M(\mathbb{C})}\\
&\le e^{3d}\Big((a_2+1)\big(\gamma(\mu)\gamma_0(x,\mu)+\gamma_1(\mu)\big)\Big),
\end{align*}
where $x\in [0,1]$ and $\mu\in P_d$.
Combining all these inequalities with
 Lemma \ref{asD} and the estimates from
 \eqref{gamma2-1}, we obtain the required representations for $D$.
\end{proof}

\begin{rem}
Note that the explicit formula for $D_0$ is the following
\[
D_{0}(x,\mu)=\left(
\begin{array}{cc}
r_1(x,\mu) & q_1(x,\mu) \\
q_2(x,\mu) & r_2(x,\mu)
\end{array}
\right),
\]
\[
q_1(x,\mu):=\int_0^x e^{i\mu (x-2t)} \sigma_1(t)\,dt,\quad
q_2(x,\mu):=\int_0^x e^{-i\mu(x-2t)} \sigma_2(t)\,dt\]
\[
r_1(x,\mu):=\int_0^x e^{i\mu (x-2t)} \tilde{\sigma}_1(x,t)\,dt,\quad
r_2(x,\mu):=\int_0^x e^{-i\mu(x-2t)} \tilde{\sigma}_2(x,t)\,dt
\]
and
$\tilde{\sigma}_j$ are given by \eqref{defSS}.
\end{rem}

\begin{rem}\label{p1}
If $p=1$, then the remainder $R_0$ from \eqref{LCdop} satisfies
\[
\|R_0(x,\mu)\|_{M(\mathbb{C})}\le
 c_2\Gamma^2(\mu),
\]
where $\Gamma$ is given by \eqref{Gam}.
\end{rem}

\section{Spectral problem}
We consider spectral problem
\begin{equation}\label{Y}
Y'(x)+J(x) Y(x)=A_\mu Y(x),\quad
 x\in [0,1],
\end{equation}
associated with the matrix problem $\eqref{v120}$
where $Y=[y_1,y_2]^T$ and
\begin{align}\label{warunki}
y_1(0)=y_2(0), \ \ y_1(1)=y_2(1).
\end{align}

Let $\mathbf{c}=\mathbf{c}(x,\mu)=[c_1,c_2]^T$ and $\mathbf{s}=\mathbf{s}(x,\mu)=[s_1,s_2]^T$ be the solutions of \eqref{Y} satisfying $c_1(0)=1$, $c_2(0)=0$ and
$s_1(0)=0$, $s_2(0)=1$. Then due to conditions \eqref{warunki} we find that the eigenvalues are the zeros of
\begin{align}\label{char}
\Phi(\lambda)=c_1(1,\lambda) + s_1(1,\lambda) - c_2(1,\lambda) - s_2(1,\lambda).
\end{align}
The eigenfunctions will be of the form:
\begin{align}
Y=[y_1,y_2]^T=[c_1(\cdot,\mu_n) + s_1(\cdot,\mu_n), c_2(\cdot,\mu_n) + s_2(\cdot,\mu_n)]^T.\label{wwcs}
\end{align}
The analysis of zeros of \eqref{char} will now lead us to characterization  of eigenvalues.

The standard approach is to derive first basic formula
for eigenvalues and then using sharp asymptotic results derive more accurate form. We thus need
results related to functions $\mathbf{s}$ and $\mathbf{c}$ from \eqref{D}.
We derive that
\begin{align}
\Phi(\mu)&=2i\sin \mu + \int_0^1 e^{(1-2t)i\mu} \Big(Q_{11}(1,t)+Q_{12}(1,t)+\sigma_1(t)\Big)dt\nonumber \\
&-\int_0^1 e^{-(1-2t)i\mu} \Big(Q_{21}(1,t)+Q_{22}(1,t)+\sigma_2(t)\Big)dt.\label{charD}
\end{align}
Changing variables in integrals we may write
\begin{align}\label{charD2}
\Phi(\mu)=2i\sin(\mu)+ V(\mu),
\end{align}
where
\begin{align}\label{Q}
V(\mu)=\int_{-1}^1e^{i\mu s}f(s)ds - \int_{-1}^1e^{-i\mu s}g(s)ds
\end{align}
and $f$, $g$ are certain function from $L_p[-1,1]$.

Note that the identities \eqref{charD} and \eqref{charD2} are true
not only for $\mu \in P_d$ but for all $\mu \in \mathbb{C}$.
It is a standard procedure (see for instance \cite{Bel}) to derive using Rouche Theorem that
zeros of $\Phi$ are in the form
$\mu_n=\pi n + \widetilde{\mu}_n$, where $(\widetilde{\mu}_n)$ is bounded. This conclusion yield that
eigenvalues lie in a certain horizontal stripe of the complex plane. We may continue and investigate more
precise the behavior of  $(\widetilde{\mu}_n)$.

The formula for $\Phi$ gives us
\begin{align}
\sin(\widetilde{\mu}_n)=\frac{(-1)^{n+1}}{2i}R(\pi n+\widetilde{\mu}_n).
\end{align}
This expression converges to zero since the convergence of the integral in \eqref{Q} follows from Lebesgue--Riemann Lemma
and the fact that $\widetilde{\mu}_n$ are bounded. Thus $\widetilde{\mu}_n\to 0$ when $n\to \infty$. Here ends the reasoning  and first claim for $p=1$.

For $1<p<2$ we may continue in order to obtain more information.
Using  $\sin x= x +\mbox{O}(x^3),$ $x\to 0$,  and the fact that $\widetilde{\mu}_n\to 0$ we obtain
\begin{align}\label{tildemu}
\widetilde{\mu}_n&=\frac{(-1)^{n+1}}{2i}\int_{-1}^1e^{i\widetilde{\mu}_n s}e^{i\pi n s}f(s)ds\\&
-\frac{(-1)^{n+1}}{2i}\int_{-1}^1e^{-i\widetilde{\mu}_n s}e^{-i\pi n s}g(s)ds +\mbox{O}(|\widetilde{\mu}_n|^3).\nonumber
\end{align}
Next, the expansion of the exponential function
$$e^{\mu t}=1+\mu t+O(|\mu|^2),\ \  \mu\to 0, \ \ |t|\le 1$$ leads to a conclusions for one of the integrals
\begin{align*}
\Big|\int_{-1}^1e^{i\widetilde{\mu}_n s}e^{i\pi n s}f(s)ds\Big|&=\int_{-1}^1e^{i\pi n s}f(s)ds\\
&+i\widetilde{\mu}_n\int_{-1}^1e^{i\pi n s}s f(s)ds+O(|\widetilde{\mu}_n|^2).
\end{align*}
Note that second integral is a product of $\widetilde{\mu}_n$ and a Fourier coefficient of the function from $L_p$, hence it would give a sequence from $l_q$, which
converges to zero. Consequently, we go back to
\eqref{tildemu} and conclude that $\widetilde{\mu}_n$ is a sum of Fourier coefficients for functions from $L_p$, hence $(\widetilde{\mu}_n)\in l_q$.
Summarizing, we showed that the eigenvalues $\mu_n$ of our spectral problem satisfy
\begin{align}\label{lqbasic}
\mu_n=\pi n + \widetilde{\mu}_n, \ \ (\widetilde{\mu}_n)\in l_q.
\end{align}

This representation for $1<p<2$ and the fact that for $p=1$ the remainder goes to zero allows us to find in both cases more accurate description of eigenvalues.
Recall we showed eigenvalues lie in $P_d$ for a certain $d>0$, thus we can use  asymptotic formulas true in a stripe.
The main tool will be the formulas for $\mathbf{c}$ and $\mathbf{s}$ and consequently for $\Phi$ from corollary \ref{Cor3.2}.

We infer that
\begin{align}
\Phi(\mu)&=2i\sin \mu + \int_0^1 e^{(1-2t)i\mu} \sigma_1(t)dt+ \int_0^1 e^{(1-2t)i\mu} \tilde{\sigma}_1(1,t)dt\nonumber \\
&-\int_0^1 e^{-(1-2t)i\mu}\sigma_2(t)dt-\int_0^1 e^{-(1-2t)i\mu}\tilde{\sigma}_2(1,t)dt+r(\mu),
\end{align}
where $$|r(\mu)|\leq c(\gamma(\mu)\gamma_0(1,\mu)+\gamma_2(\mu))\leq c(\gamma^2(\mu)+\gamma^2_0(1,\mu)+\gamma_1(\mu)).$$

The representation \eqref{lqbasic}, lemma \ref{skroty} and  discussion similar to that about
eigenvalues yield
 \begin{align}
2i(-1)^{n+1}\widetilde\mu_n &=-\int_0^1 e^{-2\pi i n t} \sigma_1(t)dt \nonumber \\
&+2(-1)^{n+1}\int_0^1\int_0^t\sigma_1(t)\sigma_2(\xi)e^{-2\pi i n t}e^{2\pi i n \xi} d \xi dt\nonumber\\
&+\int_0^1 e^{2\pi i n t}\sigma_2(t)dt+r(\mu_n),\label{charDDD}
\end{align}
For $p=1$ we have here $|r(\mu_n)|\leq c \Gamma^2(\pi n)$.

Our last aim is to prove that for $1<p<2$ there holds $(r(\mu_n))\in l_{q/2}$.
In what follows we will use a basic formula for eigenvalues \eqref{lqbasic},
a simple inequality
$|e^{iz}-1|\le |z|e^{d}$, $z\in P_d$ and the Hausdorff–Young inequality.
We infer for $\sigma\in L_p[0,1]$ that
\begin{align}
\sum_{n=1}^\infty
\Big|\int_0^xe^{\pm 2\mu_nt}\sigma(t)dt\Big|^{q}
&\le c_q\sum_{n=1}^\infty
\left|\int_0^xe^{\pm 2\pi i n t}\sigma(t)dt\right|^q\label{coeff}\\
&+
c_q\sum_{n=1}^\infty
\left(\int_0^x|e^{2i\tilde{\mu}_n t}-1||\sigma(t)|dt\right)^q \nonumber\\
&\leq c_q \|\sigma\|^q_{L_p[0,1]}+c\|\sigma\|_{L_1[0,1]}^q\sum_{n=1}^\infty
|\tilde{\mu}_n|^q\le m<\infty,\nonumber
\end{align}
for any $x\in [0,1]$.
It follows from \eqref{coeff} that
\begin{equation}\label{Nart}
\sup_{x\in [0,1]}\sum_{n=1}^\infty  \gamma_0^{q}(x,\mu_n)< \infty,
\end{equation}
Note that by \eqref{Nart}
\begin{align*}
\sum_{n=1}^\infty\gamma^q(\mu_n)
\leq c \int_0^1 \sum_{n=1}^\infty \gamma^q_0(s,\mu_n) ds < \infty.
\end{align*}
and
\begin{align*}
\|\gamma_1(\mu_n)\|_{l^{q/2}}&\le \int_0^1 |\sigma_0(s)|\|\gamma_0^2(s,\mu_n)\|_{l^{q/2}}\,ds\\
&=\int_0^1 |\sigma_0(s)| \Bigg(\sum_{n=1}^\infty \gamma_0^q(s,\mu_n)\Bigg)^{2/q}\,ds \le c\|\sigma_0\|_{L_1}.
\end{align*}
Finally,  we obtain
\begin{align*}
\sum_{n=1}^\infty |r_n|^{q/2}&< \infty.
\end{align*}

Summarizing the discussion above we proved the following fact.
\begin{thm}\label{wartif} The eigenvalues of the spectral problem \eqref{Y}-\eqref{warunki} lie in a certain stripe $P_d$ and admit the representation
\begin{align}\label{mu01}
\mu_n=\pi n + \mu_{0,n}+ \rho_n, \ \ n=1,2, \ldots
\end{align}
with
\begin{align}
\mu_{0,n} &=\frac{(-1)^{n}}{2i}\int_0^1 e^{-2\pi i n t} \sigma_1(t)dt +\frac{(-1)^{n+1}}{2i}\int_0^1 e^{2\pi i  n t}\sigma_2(t)dt\label{mu0}\\
&-i\int_0^1\int_0^t\sigma_1(t)\sigma_2(\xi)e^{-2\pi i n t}e^{2\pi i n \xi} d \xi dt\nonumber
\end{align}
and for $p=1$ there holds
\begin{align}\label{forp1}
|\rho_n|<c \Gamma^2(\pi n),
\end{align}
where $\Gamma$ is defined in \eqref{Gam},
whereas for $1<p<2$ it is true that
\begin{align}\label{forp}
\sum_{n=1}^\infty |\rho_n|^{q/2}&< \infty.
\end{align}
\end{thm}

\begin{rem}
Recall that according to lemma \ref{Jtilde} for every $x\in [0,1]$ functions $\tilde{\sigma}_j(x,\cdot)$ are from $L_r$.
If $1<p\leq\frac{4}{3}$, then $1<r\leq 2$ and Fourier coefficients of $\tilde{\sigma}_j(x,\cdot)$
are from $l_{q/2}$. Then the representation \eqref{mu01}
with
\begin{align*}
\sum_{n=1}^\infty |\rho_n|^{q/2}&< \infty
\end{align*}
 is true but with $\mu_{0,n}$ given by
\begin{align}
2i(-1)^{n+1}\mu_{0,n} &=-\int_0^1 e^{-2\pi i n t} \sigma_1(t)dt
+\int_0^1 e^{2\pi i  n t}\sigma_2(t)dt\label{mu0rem}.
\end{align}
\end{rem}

Now, we can proceed to eigenfunctions. We are going to combine results from the previous theorem with lemma \ref{asD} and
corollary \ref{Cor3.2}.

\begin{thm}\label{more}
Let $1<p<2$  and
\begin{align*}
F_1(x,t)&=\sigma_1(t)+\tilde{\sigma}_1(x,t)-(T_{\sigma_1}\tilde{\sigma}_2)(x,t)\\
F_2(x,t)&=\sigma_2(t)+\tilde{\sigma}_2(x,t)-(T_{\sigma_2}\tilde{\sigma}_1)(x,t).
\end{align*}
The eigenfunctions of the spectral problem \eqref{Y}-\eqref{warunki} admit the representation
\begin{align*}
y_1(x,\mu_n)&=e^{i\pi n x}(1+i\mu_{0,n}x)\Bigg(1+\int_0^x e^{-2\pi i n t} F_1(x,t)dt\Bigg)\\
&-2i\mu_{0,n}e^{i\pi n x}\int_0^x e^{-2\pi i n t} tF_1(x,y)dt +  r_1(x,n),\\
y_2(x,\mu_n)&=e^{-i\pi n x}(1-i\mu_{0,n}x)\Bigg(1+\int_0^x e^{2\pi i n t} F_2(x,t)dt\Bigg)\\
&+2i\mu_{0,n}e^{-i\pi n x}\int_0^x e^{2\pi i n t} tF_2(x,t)dt +  r_1(x,n),
\end{align*}
where
\begin{align*}
\sum_{n=1}^\infty \sup_{x\in[0,1]} |r_j(x,n)|^{q/2}&< \infty,
\end{align*}
\end{thm}
\begin{proof}
According to \eqref{wwcs} eigenfunctions are express by solutions $\mathbf{c}$ and $\mathbf{s}$ in the
following way
$$y_1(x,\mu_n)=c_1(x,\mu_n) + s_1(x,\mu_n)$$
and
$$y_2(x,\mu_n)=c_2(x,\mu_n) + s_2(x,\mu_n).$$
The results of lemma \ref{t1} yield that
\begin{align*}
y_1(x,\mu_n)&=e^{i\mu_n x} + \int_0^x e^{(x-2t)i\mu_n} \sigma_1(t)dt\\
&+ \int_0^x e^{(x-2t)i\mu_n} \tilde{\sigma}_1(x,t)dt-\int_0^xe^{i \mu (x-2t)}(T_{\sigma_1}\tilde{\sigma}_2)(x,t)dt+ \alpha(x,\mu_n),\\
y_2(x,\mu_n)&=e^{-i\mu_n x} + \int_0^x e^{-(x-2t)i\mu_n} \sigma_2(t)dt\\
&+ \int_0^x e^{-(x-2t)i\mu_n} \tilde{\sigma}_2(x,t)dt-\int_0^xe^{-i \mu (x-2t)}(T_{\sigma_2}\tilde{\sigma}_1)(x,t)dt+ \beta(x,\mu_n),
\end{align*}
where $$|\alpha(x,\mu_n)|+|\beta(x,\mu_n)| \leq c\gamma^2(\mu_n).$$

Repeating once more all argument used in order to derive formulas for eigenvalues, we obtain the thesis
with claimed estimates for reminders.
\end{proof}

It is possible to obtain shorter but less precise formulas for eigenfunctions. This time we use the representation \eqref{LCdop} and
comments from lemma \ref{skroty} to prove
the following fact.
\begin{cor}\label{less}
Let $1\leq p<2$, then the eigenfunctions of the spectral problem \eqref{Y}-\eqref{warunki} admit the representation
\begin{align}
y_1(x,\mu_n)&=e^{i\pi n x}\Big(1+i\mu_{0,n}x+\int_0^x e^{-2\pi i n t} \sigma_1(t)dt\nonumber\\
&+\int_0^x\int_0^s\sigma_1(s)\sigma_2(\xi)e^{-2i\mu s}e^{2i\mu \xi} d \xi ds \Big)+ r_1(x,n), \nonumber\\
y_2(x,\mu_n)&=e^{-i\pi n x}\Big(1-i\mu_{0,n}x+\int_0^x e^{2\pi i n t} (\sigma_2(t)dt\nonumber\\
&+\int_0^x\int_0^s\sigma_1(\xi)\sigma_2(s)e^{2i\mu s}e^{-2i\mu \xi} d \xi ds\Big) +  r_2(x,n),\label{wlasne}
\end{align}
where for $1<p<2$ we have
\begin{align*}
\sup_{x\in[0,1]}\sum_{n=1}^\infty  |r_j(x,n)|^{q/2}&< \infty,
\end{align*}
whereas for $p=1$ there holds
\[
|r_j(x,n)|\le
 c\Gamma^2(\pi n).
\]
\end{cor}

\section{Appendix}

\begin{lem}\label{Jtilde}
For every $x\in[0,1]$ and $j=1,2$ the functions $\tilde{\sigma}_j(x,\cdot)$ belong to $L_r[0,1]$,
therefore $\tilde{\sigma}_j\in B$, $j=1,2,$ and $\tilde{J}\in M(B)$.
\end{lem}
\begin{proof}
We take $(x,t)\in \Delta$.
Let $\widehat{\sigma}_1,\widehat{\sigma}_2$ denote the extension of $\sigma_1$ and $\sigma_2$ by zero outside $[0,1]$. Note that for every $x\in[0,1]$ we
get
\begin{align*}
\int_0^{x-t}\sigma_1(t+\xi)\sigma_2(\xi)d\xi&=\int_{-\infty}^\infty\widehat{\sigma}_1(t+\xi)\widehat{\sigma}_2(\xi)\chi(x-(t+\xi))d\xi\\
&=\int_{-\infty}^\infty\widehat{\sigma}_1(t-s)\widehat{\sigma}_2(-s)\chi(x-(t-s))ds\\
&=\Big((\widehat{\sigma}_1(\cdot)\chi(x-\cdot))\ast(\widehat{\sigma}_2(-\cdot))\Big)(t).
\end{align*}
We thus have
\begin{align*}
\Big\|\int_0^{x-t}\sigma_1(t+\xi)\sigma_2(\xi)d\xi\Big\|_{L_r[0,1]}&\leq\Big\|\int_0^{x-t}\widehat{\sigma}_1(t+\xi)\widehat{\sigma}_2(\xi)d\xi\Big\|_{L_r(\mathbb{R})}\\
&\leq \|(\widehat{\sigma}_1(\cdot)\chi(x-\cdot))\|_{L_p(\mathbb{R})}\|\widehat{\sigma}_2(-\cdot)\|_{L_p(\mathbb{R})}\\
&\leq \|\sigma_1\|_{L_p[0,1]}\|\sigma_2\|_{L_p[0,1]},
\end{align*}
hence
\begin{equation}\label{L22CY}
\|\tilde{\sigma}_j\|_{B}\le \tilde{a}, \qquad j=1,2.
\end{equation}
Clearly, a similar estimate holds  for $\tilde{\sigma}_2$ as well.

Therefore, if we consider $\epsilon$ and $x$ such that $0\leq t\leq x+\epsilon\leq 1$, then repeating the reasoning from the latter inequality, we obtain
 \begin{align*}
\Big\|&\int_0^{x+\epsilon-t}\sigma_1(t+\xi)\sigma_2(\xi)d\xi-\int_0^{x-t}\sigma_1(t+\xi)\sigma_2(\xi)d\xi\Big\|_{L_r[0,1]}\\
&\leq \|(\widehat{\sigma}_1(\cdot)\big[\chi(x+\epsilon-\cdot)-\chi(x-\cdot)\big]\|_{L_p(\mathbb{R})}\|\widehat{\sigma}_2(-\cdot)\|_{L_p(\mathbb{R})}\\
&\leq \int_{\mathbb{R}}|\widehat{\sigma}_1(s)|^p|\chi(x+\epsilon-s)-\chi(x-s)|^pds\|\widehat{\sigma}_2\|_{L_p(\mathbb{R})}.
\end{align*}
The integral in the last line converges to zero, if $\epsilon\to 0$, because of Lebesgue Theorem, hence the mapping $x\mapsto \tilde{\sigma}_j(x,\cdot)\in L_r[0,1]$ is
continuous for\linebreak  $j=1,2$.
\end{proof}

\begin{lem}\label{skroty}
The following identity holds
\begin{align*}
\int_0^x e^{-2ti\mu} \tilde{\sigma}_1(x,t)dt +
\int_0^x e^{2ti\mu} \tilde{\sigma}_2(x,t)dt&=\int_0^xe^{-2i\mu \xi}\sigma_1(\xi)d \xi \int_0^x\sigma_2(s)e^{2i\mu s} ds.
\end{align*}
Moreover, we have
\begin{align*}
\Big|\int_0^x e^{-2ti\mu} \tilde{\sigma}_1(x,t)dt +
\int_0^x e^{2ti\mu} \tilde{\sigma}_2(x,t)dt\Big|\leq c \gamma_0^2(x,\mu)
\end{align*}
and
\begin{align*}
\int_0^x e^{-2ti\mu} \tilde{\sigma}_1(x,t)dt &-
\int_0^x e^{2ti\mu} \tilde{\sigma}_2(x,t)dt\\&=-2\int_0^x\int_0^s\sigma_1(\xi)\sigma_2(s)e^{2i\mu s}e^{-2i\mu \xi} d \xi ds + \alpha_1(\mu),
\end{align*}
\begin{align*}
\int_0^x e^{2ti\mu} \tilde{\sigma}_2(x,t)dt&-\int_0^x e^{-2ti\mu} \tilde{\sigma}_1(x,t)dt\\&
=-2\int_0^x\int_0^s\sigma_1(s)\sigma_2(\xi)e^{-2i\mu s}e^{2i\mu \xi} d \xi ds + \alpha_2(\mu),
\end{align*}
where  $\alpha_j(\mu)=O(\gamma_0^2(x,\mu))$  for $j=1,2$.
\end{lem}
\begin{proof}
Note that
\begin{align*}
\int_0^x e^{-2ti\mu} \tilde{\sigma}_1(x,t)dt&=\int_0^x\int_0^s\sigma_1(s)\sigma_2(\xi)e^{-2i\mu s}e^{2i\mu \xi} d \xi ds\\
\int_0^x e^{2ti\mu} \tilde{\sigma}_2(x,t)dt&=\int_0^x\int_0^s\sigma_1(\xi)\sigma_2(s)e^{2i\mu s}e^{-2i\mu \xi} d \xi ds.
\end{align*}
Observe that the change of variables yield
\begin{align*}
\int_0^x e^{-2ti\mu} \tilde{\sigma}_1(x,t)dt&=\int_0^x\int_\xi^x\sigma_1(\xi)\sigma_2(s)e^{-2i\mu s}e^{2i\mu \xi} ds d\xi \\
&=\int_0^x\int_s^x\sigma_1(s)\sigma_2(\xi)e^{-2i\mu \xi}e^{2i\mu s} d\xi ds,
\end{align*}
thus
\begin{align*}
\int_0^x e^{-2ti\mu} \tilde{\sigma}_1(x,t)dt &+
\int_0^x e^{2ti\mu} \tilde{\sigma}_2(x,t)dt\\&=\int_0^x\int_0^x\sigma_1(\xi)\sigma_2(s)e^{2i\mu s}e^{-2i\mu \xi} d \xi ds\\
&=\int_0^xe^{-2i\mu \xi}\sigma_1(\xi)d \xi \int_0^x\sigma_2(s)e^{2i\mu s} ds.
\end{align*}
This step shows that
\begin{align*}
\Big|\int_0^x e^{-2ti\mu} \tilde{\sigma}_1(x,t)dt +
\int_0^x e^{2ti\mu} \tilde{\sigma}_2(x,t)dt\Big|\leq c \gamma_0^2(x,\mu).
\end{align*}

What is more,
then
\begin{align*}
\int_0^x e^{-2ti\mu} \tilde{\sigma}_1(x,t)dt &-
\int_0^x e^{2ti\mu} \tilde{\sigma}_2(x,t)dt\\&=\int_0^xe^{-2i\mu \xi}\sigma_1(\xi)d \xi \int_0^x\sigma_2(s)e^{2i\mu s} ds\\
&-2\int_0^x\int_0^s\sigma_1(\xi)\sigma_2(s)e^{2i\mu s}e^{-2i\mu \xi} d \xi ds,
\end{align*}
thus
\begin{align*}
\int_0^x e^{-2ti\mu} \tilde{\sigma}_1(x,t)dt &-
\int_0^x e^{2ti\mu} \tilde{\sigma}_2(x,t)dt\\&=-2\int_0^x\int_0^s\sigma_1(\xi)\sigma_2(s)e^{2i\mu s}e^{-2i\mu \xi} d \xi ds\\
& + \alpha_1(\mu),
\end{align*}
where $\alpha_1(\mu)=O(\gamma_0^2(x,\mu))$.
Analogously we get  the last claim
\end{proof}

\begin{lem}\label{Tsigma}
The linear
operator $T_{\sigma}$
\begin{align}\label{TT}
(T_{\sigma}f)(x,t)=\int_0^{x-t}\sigma(t+\xi)f(t+\xi,\xi)d\xi=
\int_t^x\sigma(s)f(s,s-t)ds, 
\end{align}
where $\sigma\in L_p[0,1]$is  bounded in $B$.
\begin{proof}
Note that
\begin{align}
\Big(\int_0^x|(T_\sigma f)(x,t)|^rdt\Big)^{1/r}&=\Big(\int_0^x\Big|\int_0^x\chi(s-t)\sigma(s)f(s,s-t)ds\Big|^rdt\Big)^{1/r}\nonumber\\
&\leq \int_0^x|\sigma(s)|\Big(\int_0^s|f(s,s-t)|^rdt\Big)^{1/r}ds\nonumber \\
&\leq \int_0^x|\sigma(s)|\Big(\int_0^s|f(s,\tau)|^rd\tau\Big)^{1/r}ds\nonumber \\
&\leq \int_0^x|\sigma(s)|ds \sup_{s\in[0,1]}\Big(\int_0^s|f(s,\tau)|^rd\tau\Big)^{1/r}\nonumber \\
&\leq \|\sigma\|_{L_1}\|f\|_B.\label{Ts}
\end{align}

For the proof of continuity we take  $\epsilon$ and $x$ such that $0\leq t\leq x+\epsilon\leq 1$. Then
\begin{align*}
\Big\|(T_\sigma f)(x+\epsilon,\cdot)-(T_\sigma f)(x,\cdot)\Big\|_{L_r[0,1]}&\leq\Big(\int_0^x\Big|\int_x^{x+\epsilon}\sigma(s)f(s,s-t)ds\Big|^rdt\Big)^{1/r}\nonumber\\
&+\Big(\int_x^{x+\epsilon}\Big|\int_t^{x+\epsilon}\sigma(s)f(s,s-t)ds\Big|^rdt\Big)^{1/r}.
\end{align*}
First integral may be estimated as follows
\begin{align*}
\Big(\int_0^x\Big|\int_x^{x+\epsilon}\sigma(s)f(s,s-t)ds\Big|^rdt\Big)^{1/r}
&\leq \int_x^{x+\epsilon}|\sigma(s)|\Big(\int_0^x|f(s,s-t)|^rdt\Big)^{1/r}ds\\
&\leq \int_x^{x+\epsilon}|\sigma(s)|\Big(\int_0^s|f(s,\tau)|^rd\tau\Big)^{1/r}ds\\
&\leq\|f\|_B \int_x^{x+\epsilon}|\sigma(s)|ds
\end{align*}
and this expression goes to zero whenever $\epsilon$ does.

Second integral can be treated in  an analogous way, hence the proof is completed
\end{proof}
\end{lem}

\begin{lem}\label{lemT}
The operators $T_{kj}, k,j=1,2,\quad k\not=j$ satisfy the following estimate
\begin{align}\label{Tn}
\|T_{kj}^nf\|_{B}
\le
\frac{a^n}
{n!}\|f\|_{B},\qquad f\in B,\quad n\in \mathbb{N},\quad k,j=1,2,\quad k\not=j.
\end{align}
\end{lem}

\begin{proof} Consider the operator $T_{12}$.
Note that directly from third line of \eqref{Ts} we get
\begin{align}
\Big(\int_0^x|(T_{12}f)(x,t)|^rdt\Big)^{1/r}
&\leq \int_0^x|\sigma_1(s)|\Big(\int_0^s|(T_{\sigma_2}f)(s,\tau)|^rd\tau\Big)^{1/r}ds\nonumber\\
& \leq  \int_0^x|\sigma_1(s)|\int_0^s|\sigma_2(\tau)| \Big(\int_0^\tau|f(\tau,\xi)|^rd\xi\Big)^{1/r} d\tau ds\nonumber\\
& \leq  \|f\|_B\int_0^x|\sigma(s)|\int_0^s|\sigma_2(\tau)|ds\label{Tss}.
\end{align}
Define $\eta\in C[0,1]$ by
\[
\eta(x):=\int_0^x |\sigma_1(s)|\left(\int_0^s |\sigma_2(\tau)|\,d\tau\right)\,ds,\qquad x\in [0,1].
\]
This function is increasing and bounded by $a=\|\sigma_1\|_{L_1}\|\sigma_2\|_{L_1}$. It suffices to prove that
for all $(x,t)\in \Delta$ and $n=1,2,\dots,$
\begin{equation}\label{100}
\Big(\int_0^x|(T^n_{12}f)(x,t)|^rdt\Big)^{1/r}
\le \frac{\|f\|_{B}}{n!}
\eta^n(x),\qquad
 f\in B.
\end{equation}
For $n=1$ the estimate \eqref{100}
was shown above.
Arguing by induction, suppose that \eqref{100} holds for some $n\in \mathbb{N}$.
Then for $(x,t)\in \Delta$ and $f\in B$ from \eqref{Tss} we have
\begin{align*}
\Big(\int_0^x|(T^{n+1}_{12}f)&(x,t)|^rdt\Big)^{1/r}\\ &\leq  \int_0^x|\sigma_1(s)|\int_0^s|\sigma_2(\tau)| \Big(\int_0^\tau|(T^{n}_{12}f)(\tau,\xi)|^rd\xi\Big)^{1/r} d\tau ds\nonumber\\
&\le \frac{\|f\|_{B}}{n!}
\int_0^x |\sigma_1(s)| \int_{0}^s |\sigma_2(\tau)| \eta^n(\tau)\,d\tau\,ds
\\
&\le \frac{\|f\|_{B}}{n!}
\int_0^x |\sigma_1(s)| \int_0^s |\sigma_2(\tau)|\,d\tau\,\eta^n (s)\,ds
\\
&=\frac{\|f\|_{B}}{n!}
\int_0^x \eta^n(s)\,d\eta(s)
=\frac{\|f\|_{B}}{(n+1)!}\eta^{n+1}(x).
\end{align*}
Therefore \eqref{100} hold true and then after taking supremum over $x\in[0,1]$ we get \eqref{Tn}.
\end{proof}

Next proposition we state below without a proof, since it can be found in \cite[Prop. 6.1]{GRZ1}.
\begin{prop}\label{4.9}
If $\sigma_j\in L_p[0,1]$, $1\leq p <2 $ and $F\in M(B)$, then
\begin{equation}\label{ind0}
\int_0^xe^{-2i\mu t}(\tilde{T} F)(x,t)dt
=-\int_0^xe^{-2i\mu s}J(s)\int_0^{s}e^{2i\mu \xi}F(s,\xi)d\xi ds.
\end{equation}
Moreover,
\begin{align}\label{wyr2}
\int_0^xe^{-2i\mu t} &(\tilde{T}\tilde{J})(x,t)\,dt \nonumber\\
&=-\int_0^xe^{2i\mu y} \Bigg(\int_y^x J(z)e^{-2i\mu z}\, d z\,
\int_0^{y}J^T(\tau)e^{-2i\mu \tau} d\tau
\Bigg)J^T(y) dy.
\end{align}
\end{prop}

\begin{lem}\label{Asimp}
If $\mu\in P_d,$ then there hold the following inequalities
\begin{equation}\label{EstIm0}
\left\|\int_0^x e^{-2tA_\mu}\tilde{J}(x,t)\,dt\right\|_{M(C[0,1])}
\le 2e^{2d}\tilde{a}_0\,\gamma(\mu),
\end{equation}
\begin{equation}\label{EstIm100}
\left\|\int_0^x e^{-2tA_\mu}(\tilde{T}\tilde{J})(x,t)\,dt
\right\|_{M(C[0,1])}
\le 2 e^{4d} \tilde{a}_0^2\gamma(\mu),
\end{equation}
\begin{align}\label{EstIm1}
\left\|\int_0^x e^{-2tA_\mu}(\tilde{T}\tilde{J})(x,t)\,dt\right\|_{M( \mathbb{C})}
&\le 2(a_2+1)e^{2d}
\Big(\gamma(\mu)\gamma_0(x,\mu)\nonumber+\gamma_1(\mu)\Big),\\&
 x\in [0,1],
\end{align}
\begin{equation}\label{EstIm3}
\left\|\int_0^x e^{-2tA_\mu}(\tilde{T}^n \tilde{J})(x,t)\,dt
\right\|_{M(C[0,1])} \le 2e^{2n d}
\frac{a_1^{n-2}}{(n-2)!}
\gamma_2(\mu),\quad n\ge 2.
\end{equation}
\end{lem}

\begin{proof}
Note that
\begin{align}
\left\|\int_0^xe^{-2i\mu t}
\tilde{J}(x,t)dt\right\|_{M(\mathbb{C})}&=
\left\|\int_0^x J(s)e^{-2i\mu s}
\int_0^s J(\xi)e^{2i\mu \xi} d \xi ds\right\|_{M( \mathbb{C})}\nonumber
\\
&=\left|\int_0^x e^{-2i\mu s}\sigma_1(s)
\int_0^s e^{2i\mu \xi}\sigma_2(\xi) d \xi ds\right|\nonumber \\&+
\left|\int_0^x e^{-2i\mu s}\sigma_2(s)
\int_0^s e^{2i\mu \xi}\sigma_1(\xi) d \xi ds\right|
\label{EstIm00}\\
&\le  e^{2d}\Big\{\|\sigma_1\|_{L_p}\Big\|\int_0^s e^{2i\mu \xi}\sigma_2(\xi)d\xi \Big\|_{L_q}
\nonumber\\&+\|\sigma_2\|_{L_p}\Big\|\int_0^s e^{2i\mu \xi}\sigma_1(\xi)d\xi \Big\|_{L_q}\Big\}\nonumber
\\
&\le  e^{2d}\max\{\|\sigma_1\|_{L_p},\|\sigma_2\|_{L_p}\}\,\gamma(\mu),\quad x\in [0,1].\nonumber
\end{align}
We thus proved the estimate \eqref{EstIm0}.

Next,
from \eqref{ind0}, if $\mu\in P_d$, $x\in [0,1]$ and $F\in M(B)$, then
\begin{equation}\label{Nes0}
\left\|\int_0^xe^{-2i\mu t}(\tilde{T} F)(x,t)dt\right\|_{M( \mathbb{C})}
\le e^{2d}\int_0^x \left\|J(s)
\int_0^{s}e^{2i\mu \xi}F(s,\xi)d\xi\right\|_{M(\mathbb{C})}\,ds,
\end{equation}
and
\begin{equation}\label{Nes1}
\left\|\int_0^xe^{-2i\mu t}(\tilde{T} F)(x,t)dt\right\|_{M(C[0,1])}
\le e^{2d}a_0
\left\|\int_0^{s}e^{2i\mu \xi}F(s,\xi)d\xi\right\|_{M(C[0,1])}.
\end{equation}
We use \eqref{Nes1} and \eqref{EstIm00} to obtain that
\begin{align*}
\left\|\int_0^x e^{-2i\mu t}(\tilde{T}\tilde{J})(x,t)\,dt
\right\|_{M(C[0,1])}
&\le e^{2d}a_0
\left\|\int_0^{s}e^{2i\mu \xi}\tilde{J}(s,\xi)\,d\xi\right\|_{M(C[0,1])}\\
&\le e^{4d}\tilde{a}_0^2\gamma(\mu),
\end{align*}
thus, the estimate
\eqref{EstIm100} holds.

Due to the estimate
\[
\left|\int_0^x \sigma_0(s)
\gamma_0(y,\mu)\,dy\right|\le \|\sigma_0\|_{L_p}\left\|\gamma_0(y,\mu)\right\|_{L_q}
\le a_2\gamma(\mu),
\]
the inequality \eqref{EstIm1} holds if
\begin{equation}\label{EstIm10}
\Bigg\|
\int_0^xe^{-2i\mu t} (\tilde{T}\tilde{J})(x,t)\,dt\Bigg\|_{M(\mathbb{C})}
\le
 e^{2d}\left(
\gamma_1(\mu)
+\gamma_0(x,\mu)
\int_0^x \sigma_0(s)
\gamma_0(y,\mu)\,dy\right).
\end{equation}
Whereas using \eqref{wyr2}, \eqref{EstIm0}, we have
\begin{align*}
\Bigg\|
\int_0^x&e^{-2i\mu t} (\tilde{T}\tilde{J})(x,t)\,dt\Bigg\|_{M({\rm C})}\\
&\le
e^{2d}
\int_0^x\left\|\int_{y}^x e^{-2i\mu z}J(z) d z \int_0^y e^{-2i\mu \tau}J^T(\tau)\,d\tau J^T(y)\right\|_{M_2({\rm C})}
dy
\\
&\le e^{2d}
\int_0^x \sigma_0(y)
\left\|\int_{y}^x e^{-2i\mu z}J(z) d z\right\|_{M({\rm C})}
\left\| \int_0^y e^{-2i\mu \tau}J(\tau)\,d\tau \right\|_{M({\rm C})}
dy
\\
&\le e^{2d}
\int_0^x \sigma_0(y)
\left\| \int_0^y e^{-2i\mu \tau}J(\tau)\,d\tau \right\|^2_{M({\rm C})}
dy
\\
&+e^{2d}\left\|\int_0^x e^{-2i\mu z}J(z) d z\right\|_{M({\rm C})}
\int_0^x \sigma_0(s)
\left\| \int_0^y e^{-2i\mu \tau}J(\tau)\,d\tau \right\|_{M({\rm C})}
dy
\\
&\le e^{2d}
\gamma_1(\mu)
+e^{2d}\gamma_0(x,\mu)
\int_0^x \sigma_0(s)
\gamma_0(y,\mu)\,dy,
\end{align*}
and \eqref{EstIm10} follows.

The estimate \eqref{EstIm3} will be showed, if we prove that
for all $n\ge 2$ and any $x\in [0,1],$
\begin{align}\label{EstIm30}
\left\|\int_0^x e^{-2i\mu t}(\tilde{T}^n \tilde{J})(x,t)\,dt\right\|_{M({\rm C})}
\le \frac{e^{2n d}}{(n-2)!}
\left(\int_0^x \sigma_0(s)\,ds\right)^{n-2}
\gamma_2(\mu).
\end{align}
We proceed  by induction.
Using \eqref{Nes0}
for $F=\tilde{T}\tilde{J}$
and  \eqref{EstIm10},
we note that
\begin{align*}
 \Big\|\int_0^x e^{-2i\mu t}(\tilde{T}^2\tilde{J})&(x,t)\,dt\Big\|_{M({\rm C})}\\&
\le e^{2d}\int_0^x \sigma_0(s)
\left\|\int_0^{s}e^{2i\mu \xi}(\tilde{T}\tilde{J})(s,\xi)d\xi\right\|_{M({\rm C})}\, ds
\\
&\le e^{4d}\int_0^x \sigma_0(s)
\left(\gamma_0(s,\mu)\int_0^s\sigma_0(y)
\gamma_0(y,\mu)\,dy+\gamma_1(\mu)\right)\, ds
\\
&\le e^{4d}\int_0^x \sigma_0(s)\gamma_0(s,\mu)
\int_0^{s}\sigma_0(y)
\gamma_0(y,\mu)\,dy\, ds
+e^{4d}a_1\gamma_1(\mu)
\\
&\le
 e^{4d}\frac{\left(\int_0^x \sigma_0(s)\gamma_0(s,\mu)\,ds\right)^2}{2}
+e^{4d}a_1\gamma_1(\mu)
\\
&\le e^{4d}\big(a^2_2\gamma^2(\mu)+e^{4d}a_1\gamma_1(\mu)\big).
\end{align*}
Therefore, \eqref{EstIm30} holds for $n=2$.

Let suppose now that \eqref{EstIm30}
holds for some $n\ge 2$.
We thus once again use \eqref{Nes0} to derive
\begin{align*}
\Big\|\int_0^x e^{-2i\mu t}(\tilde{T}^{n+1}\tilde{J})&(x,t)\,dt\Big\|_{M({\rm C})}\\
&\le e^{2d} \int_0^x \sigma_0(s)
\left\|\int_0^{s}e^{2i\mu \xi}(\tilde{T}^n\tilde{J})(s,\xi)d\xi
\right\|_{M({\rm C})} ds
\\
&\le \frac{e^{2(n+1) d}}{(n-2)!}\gamma_2(\mu)\int_0^x\sigma_0(s)
 \left(\int_0^s\sigma_0(\tau)\,d\tau\right)^{n-2}\,ds
\\
&=\frac{e^{2(n+1)d}\gamma_2(\mu)}{(n-1)!}
 \left(\int_0^x \sigma_0(\tau)\,d\tau\right)^{n-1},\quad x\in [0,1],
\end{align*}
thus \eqref{EstIm30}
holds also for  $n+1,$ and  the proof
of \eqref{EstIm30} is completed.
\end{proof}

\bibliographystyle{plain}
\bibliography{bibi2}

\end{document}